\documentclass[notitlepage,a4paper] {article}

\usepackage [english] {babel} 
\usepackage [latin1]{inputenc}

\usepackage {amsmath}
\usepackage {amsthm}
\usepackage {amssymb}
\usepackage {amsfonts}
\usepackage {mathrsfs}

\usepackage{url}
\usepackage {enumerate}

\usepackage{tikz}
\usepackage{graphicx}

\usetikzlibrary{matrix,positioning}

\usepackage{float}
\usepackage{caption}

\theoremstyle{plain}
\newtheorem*{mainTeorema*}{ {\bf Main Theorem} }
\newtheorem{teorema}{ {\bf Theorem} }[section]
\newtheorem{lema}[teorema]{ {\bf Lemma} }

\theoremstyle{definition}
\newtheorem{definicion}{ {\bf Definition} }[section]

\theoremstyle{remark}

\newtheorem{nota}{ {\bf Remark} }[section]

\newenvironment{acknowledgements}{
{ \flushleft{ \large{ \textbf{ \textrm{ Acknowledgements } } } } }
  \newline
  \newline 
}{}

\addto\captionsenglish{}

\title{\large \bf THE BINARY GOLDBACH CONJECTURE IS ALSO TRUE}
\author{ {\sl Ricardo G. Barca}}
\date{}

\begin{document} 

\maketitle

\begin{abstract}

The binary Goldbach conjecture asserts that every even integer greater than $4$ is the sum of two primes. In a preceding paper we have proved that there exists a positive integer $K_\alpha$ such that every even integer $x > p_k^2$ can be expressed as the sum of two primes, where $p_k$ is the $k$th prime number and $k > K_\alpha$. In this paper we provide an estimation for $K_\alpha$, and from this result it follows that the binary Goldbach conjecture is true.

\end{abstract}

% seccion 1
%-------------------------------------------------------------------------------------------

\section{Introduction} \label{seccion01}

This paper is the continuation of the author's work~\cite{rgb} and proves the binary Goldbach conjecture.

\begin{mainTeorema*}

Every even integer $x$ greater than $4$ can be expressed as the sum of two primes.

\end{mainTeorema*}

Let $S_k \; (k \ge 4)$ be a given partial sum of the series $\sum s_k$. Recall the notation $m_k$ for the period of the partial sum $S_k$. For every partial sum $S_h$ from level $h=1$ to level $h=k$, let us consider the interval $I[1, m_k]_h$, the Left interval $I[1, p_k^2]_h$ and the Right interval $I[p_k^2+1, m_k]_h$. Furthermore, recall the notation $\delta_h$, $\delta_h^{ L_k }$ and $\delta_h^{ R_k }$ for the density of permitted $h$-tuples within the intervals $I[1, m_k]_h$, $I[1, p_k^2]_h$ and $I[p_k^2+1, m_k]_h$, respectively.

In the previous paper~\cite{rgb} we showed that there exists $K_\alpha > 4$ such that every even number greater than $p_k^2 \; (k > K_\alpha)$ is the sum of two primes. In the present paper we shall provide an estimation of $K_\alpha$, proving this way the Goldbach conjecture.

We need to complement the discussion in~\cite{rgb} about the behaviour of $\delta_h^{ L_k }$ and $\delta_h^{ R_k }$ (in the Left and Right blocks of the partition of the first period of $S_k$, respectively) with additional observations, that we shall make in the following subsections.

\subsection{\small The behaviour of $\delta_h^{ R_k }$ in the Right block of the partition.} \label{subsec11}

Recall the notation $c_h^{ L_k }$ to denote the number of permitted $h$-tuples within the Left interval $I[1, p_k^2]_h$, and the notation $c_h^{ R_k }$ to denote the number of permitted $h$-tuples within the Right interval $I[p_k^2+1, m_k]_h$; furthermore, recall the notation $c_h$ for the number of permitted $h$-tuples within a period of the partial sum $S_h$ and $c_h^{\prime}$ to denote the number of permitted $h$-tuples within the interval $I [1, m_k ]_h$ of every partial sum $S_h \; (1 \le h \le k)$. We denote by $\overline{ c_h^{L_k} }$ and $\overline{ c_h^{R_k} }$ the average number of permitted $h$-tuples within the intervals $I[1, p_k^2]_h$ and $I[p_k^2+1, m_k]_h \; (1 \le h \le k)$, respectively.

We begin by the following lemma.

% Lema 1.1

\begin{lema} \label{S1pk2o(ck)}

Let $S_k$ be a given partial sum. Let $p_k$ be the characteristic prime modulus of the partial sum $S_k$. Let $c_k$ be the number of permitted $k$-tuples within the period of $S_k$. We have $p_k^2 = {\rm o} ( c_k )$.

\begin{proof}

Using~\cite[Proposition 2.3]{rgb}, we have

\begin{align} \label{S1EQpk2/ck}
&\frac{p_k^2}{c_k} = \frac{p_k^2}{\left ( p_1 - 1 \right ) \left ( p_2 - 2 \right ) \left ( p_3 - 2 \right ) \cdots \left ( p_k - 2 \right )} = \\ \nonumber &= \left ( \frac{1}{\left ( p_1 - 1 \right ) \left ( p_2 - 2 \right ) \left ( p_3 - 2 \right ) \cdots \left ( p_{k-2} - 2 \right )} \right ) \left ( \frac{p_k}{p_{k-1} - 2} \right ) \left ( \frac{p_k}{p_k - 2} \right ).
\end{align}

Let $g_{k-1}$ denote the gap $p_k - p_{k-1}$; so, $p_k / ( p_{k-1} - 2 ) = ( p_{k-1} + g_{k-1} ) / ( p_{k-1} - 2 )$. By the Bertrand--Chebyshev theorem, we have $g_{k-1} <  p_{k-1} \Longrightarrow ( p_{k-1} + g_{k-1} ) / ( p_{k-1} - 2 ) < 2 p_{k-1} / ( p_{k-1} - 2 )$. It follows that $\lim_{k \to \infty} p_k / ( p_{k-1} - 2 ) < \lim_{k \to \infty} 2 p_{k-1} / ( p_{k-1} - 2 ) = 2$. Since $\lim_{k \to \infty} p_k / ( p_k - 2 ) = 1$, returning to \eqref{S1EQpk2/ck}, clearly, $\lim_{k \to \infty} p_k^2 / c_k = 0$.

\end{proof}
\end{lema}

Since $p_k^2 = {\rm o } ( m_k )$, for every partial sum $S_h$ from $h=1$ to $h=k$, the size of the Right interval $I[p_k^2+1, m_k]_h$ approximates the size of the interval $I[1, m_k]_h$ more and more closely as the level $k$ increases. On the other hand $p_k^2 = {\rm o} ( c_k )$, by Lemma~\ref{S1pk2o(ck)}; furthermore, $c_h^\prime > c_k = c_k^\prime \; (1 \le h < k)$ by~\cite[Proposition 2.3]{rgb} and~\cite[Lemma 6.1]{rgb}. Therefore, the proportion of permitted $h$-tuples within every Right interval $I[p_k^2+1, m_k]_h \; (1 \le h \le k)$ approximates the proportion of permitted $h$-tuples within the respective interval $I[1, m_k]_h$ more and more closely as the level $k$ increases. Therefore, as the level $k$ increases, we expect that every value of $\delta_h^{ R_k } \; (1 \le h \le k)$ approximates the respective average $\delta_h$ more and more closely, regardless of the combination of selected remainders in the sequences $s_h$ that form the partial sum $S_k$.

The following lemma shows that as $k\to \infty$, the $h$-density within the Right interval $I[p_k^2+1, m_k]_h$ of every partial sum $S_h \; (1 \le h \le k)$ converges uniformly to the average $\delta_h$.

% Lema 1.2

\begin{lema} \label{S1deltahRkTodeltah}

Let $S_k \; (k \ge 4)$ be a partial sum of the series $\sum s_k$. Let us consider the Right interval $I[p_k^2+1, m_k]_h$ in every partial sum $S_h$ from level $h = 1$ to level $h = k$. For every $\epsilon > 0$, there exists $N$ (depending only on $\epsilon$) such that $k > N$ implies $| \delta_h^{ R_k } - \delta_h | < \epsilon$, for every partial sum $S_h$ from level $h = 1$ to level $h = k$, regardless of the combination of selected remainders in the sequences $s_h$ that form every partial sum $S_k$.

\begin{proof}

\begin{enumerate}[Step 1.]

\item The size of the Right interval $I[p_k^2+1, m_k]_h$ of the partial sum $S_h$, by definition, is equal to $m_k - p_k^2$, so the number of subintervals of size $p_h$ within the Right interval is equal to $(m_k - p_k^2) / p_h \; (1 \le h \le k)$. Denoting by $c_h^{ R_k}$ the number of permitted $h$-tuples within $I[p_k^2+1, m_k]_h$, by definition, we have

\begin{align} \label{S1EQdeltaRh}
\delta_h^{ R_k } = \frac{ c_h^{ R_k } }{ \left (m_k - p_k^2 \right ) / p_h } \qquad (1 \le h \le k).
\end{align}

\item Let us denote by $m_h$ the period of the partial sum $S_h$ and by $c_h$ the number of permitted $h$-tuples within a period of the partial sum $S_h$. For every level from $h = 1$ to $h = k$, let $c_h^\prime$ be the number of permitted $h$-tuples within the interval $I[1, m_k]_h$ of the partial sum $S_h$. Using~\cite[Lemma 6.1]{rgb}, we obtain

\begin{align} \label{S1EQchprima1}
&c_1^\prime = c_1 p_2 p_3 \cdots p_k, \\ \nonumber
&c_2^\prime = c_2 p_3 p_4 \cdots p_k, \\ \nonumber
&\ldots \\ \nonumber
&c_h^\prime = c_h p_{h+1} p_{h+2} \cdots p_k, \\ \nonumber
&\ldots \\ \nonumber
&c_k^\prime = c_k.
\end{align}

Note that $c_h^\prime$ increases as the level decreases from $h = k$ to $h = 1$ (see~\cite[Proposition 2.3]{rgb}). For every level from $h = 1$ to $h = k$, since $I[1, m_k]_h = I[ 1, p_k^2]_h \cup I[p_k^2+1, m_k]_h$, the number of permitted $h$-tuples within the Right interval $I[p_k^2+1, m_k]_h$ cannot exceed $c_h^\prime$, so we have $c_h^{ R_k } \le c_h^\prime$. On the other hand, the number of permitted $h$-tuples within the Left interval $I[ 1, p_k^2]_h$ of the partial sum $S_h$ cannot exceed the size $p_k^2$ of the Left interval. Therefore, $c_h^\prime - p_k^2 \le c_h^{ R_k }$. Consequently, replacing the numerator in \eqref{S1EQdeltaRh} by $c_h^\prime - p_k^2$ and $c_h^\prime$, we obtain

\begin{align*}
\frac{ c_h^\prime - p_k^2 }{ \left (m_k - p_k^2 \right ) / p_h } \le \delta_h^{ R_k } \le \frac{ c_h^\prime }{ \left (m_k - p_k^2 \right ) / p_h } \qquad (1 \le h \le k),
\end{align*}

so, extracting the common factors $c_h^\prime$ and $m_k$, we have

\begin{align*}
\frac{ c_h^\prime }{ m_k / p_h } \left ( \frac{ 1 - p_k^2 / c_h^\prime }{ 1 - p_k^2 / m_k } \right ) \le \delta_h^{ R_k } \le \frac{ c_h^\prime }{ m_k / p_h } \left ( \frac{ 1 }{ 1 - p_k^2 / m_k } \right ).
\end{align*}

Now, by definition, $m_k = p_1 p_2 p_3 \cdots p_h p_{h+1} p_{h+2} \cdots p_k = m_h p_{h+1} p_{h+2} \cdots p_k$. Then, using \eqref{S1EQchprima1} and canceling common factors, we obtain

\begin{align*}
\frac{ c_h }{ m_h / p_h } \left ( \frac{ 1 - p_k^2 / c_h^\prime }{ 1 - p_k^2 / m_k } \right ) \le \delta_h^{ R_k } \le \frac{ c_h }{ m_h / p_h } \left ( \frac{ 1 }{ 1 - p_k^2 / m_k } \right ).
\end{align*}

By definition,

\begin{align*}
\delta_h = \frac{ c_h }{ m_h / p_h }.
\end{align*}

Therefore, for every partial sum $S_h$ from level  $h = 1$ to level $h = k$, regardless of the combination of selected remainders in the sequences $s_h$ that form the partial sum $S_k$, we have the bounds

\begin{align} \label{S1EQboundsDeltaRh}
\delta_h \left ( \frac{ 1 - p_k^2 / c_h^\prime }{ 1 - p_k^2 / m_k } \right ) \le \delta_h^{ R_k } \le \delta_h \left (\frac{ 1 }{ 1 - p_k^2 / m_k } \right ).
\end{align}

\item Now, let $\epsilon > 0$ be a given small number, and let $N \ge 12$. For level $h=k$, from \eqref{S1EQboundsDeltaRh} we have

\begin{align*}
\delta_k \left ( \frac{ 1 - p_k^2 / c_k^\prime }{ 1 - p_k^2 / m_k } \right ) \le \delta_k^{ R_k } \le \delta_k \left (\frac{ 1 }{ 1 - p_k^2 / m_k } \right ).
\end{align*}

On the one hand, $p_k^2 = {\rm o}( m_k )$; on the other hand, $c_k = c_k^\prime$, by \eqref{S1EQchprima1}, thus, by Lemma~\ref{S1pk2o(ck)}, $p_k^2 = {\rm o} ( c_k^\prime )$. Moreover, it follows from~\cite[Proposition 2.3]{rgb} that $c_k^\prime < m_k$. Therefore, we can take sufficiently large $N$ that for level $k > N$,

\begin{align} \label{S1EQinequalitiesRk}
\delta_k - \frac{\epsilon}{2} < \delta_k \left ( \frac{ 1 - p_k^2 / c_k^\prime }{ 1 - p_k^2 / m_k } \right ) \le \delta_k^{ R_k } \le \delta_k \left (\frac{ 1 }{ 1 - p_k^2 / m_k } \right ) < \delta_k + \frac{\epsilon}{2},
\end{align}

at level $h=k$.

\item Now, the rightmost inequality in \eqref{S1EQinequalitiesRk} implies

\begin{align*}
\delta_k \left ( \frac{ 1 }{ 1 - p_k^2 / m_k } - 1 \right ) < \frac{\epsilon}{2}.
\end{align*}

For a given level $h<k$, since $k > N \ge 12$ by assumption, it is easy to check using~\cite[Lemma 3.2]{rgb} and~\cite[Corollary 3.3]{rgb} that $\delta_h \le \delta_k$. Hence,

\begin{align} \label{S1EQrightInequality}
\delta_h \left ( \frac{ 1 }{ 1 - p_k^2 / m_k } - 1 \right ) < \frac{\epsilon}{2} \Longrightarrow \delta_h \left (\frac{ 1 }{ 1 - p_k^2 / m_k } \right ) < \delta_h + \frac{\epsilon}{2}.
\end{align}

\item The leftmost inequality in \eqref{S1EQinequalitiesRk} implies

\begin{align*}
\delta_k \left ( 1 - \frac{ 1 - p_k^2 / c_k^\prime }{ 1 - p_k^2 / m_k } \right ) < \frac{\epsilon}{2}.
\end{align*}

For a given level $h<k$, since $k > N \ge 12$, we have $\delta_h \le \delta_k$ (see Step 4). On the other hand, $c_k^\prime = c_k < c_h^\prime < m_k$, where $1 \le h < k$ (see \eqref{S1EQchprima1} and~\cite[Proposition 2.3]{rgb}). Hence, replacing $\delta_k$ with $\delta_h$ and $c_k^\prime$ with $c_h^\prime$, we obtain

\begin{align} \label{S1EQleftInequality}
\delta_h \left ( 1 - \frac{ 1 - p_k^2 / c_h^\prime }{ 1 - p_k^2 / m_k } \right ) < \frac{\epsilon}{2} \Longrightarrow \delta_h - \frac{\epsilon}{2} < \delta_h \left ( \frac{ 1 - p_k^2 / c_h^\prime }{ 1 - p_k^2 / m_k } \right ).
\end{align}

\item We now prove the lemma. By \eqref{S1EQboundsDeltaRh}, \eqref{S1EQinequalitiesRk}, \eqref{S1EQrightInequality}, and \eqref{S1EQleftInequality}, for $k > N$, we can write

\begin{align*}
\delta_h - \frac{\epsilon}{2} < \delta_h \left ( \frac{ 1 - p_k^2 / c_h^\prime }{ 1 - p_k^2 / m_k } \right ) \le \delta_h^{ R_k } \le \delta_h \left (\frac{ 1 }{ 1 - p_k^2 / m_k } \right ) < \delta_h + \frac{\epsilon}{2},
\end{align*}

for every level from $h = 1$ to $h = k$. This result implies $| \delta_h^{ R_k } - \delta_h | < \epsilon$ for every level from $h = 1$ to $h = k \; (k > N)$, regardless of the combination of selected remainders in the sequences $s_h$ that form every partial sum $S_k$.

\end{enumerate}

\end{proof}
\end{lema}

\subsection{\small The behaviour of $\delta_h^{ L_k }$ in the Left block of the partition.} \label{subsec12}

Let us consider every partial sum $S_h$ from $h=1$ to $h=k$. Although the increase in the number of permitted $h$-tuples within one interval is equal to the decrease in the number of permitted $h$-tuples within the other interval, the increase in the $h$-density within one interval is not equal to the decrease in the $h$-density within the other interval. The following lemma gives the relationship between the $h$-density within $I[1, p_k^2]_h$ (denoted by $\delta_h^{ L_k }$) and the $h$-density within $I[p_k^2+1, m_k]_h$ (denoted by $\delta_h^{ R_k }$). Recall the notation $\{ \delta_h^{ L_k } \}$ to denote the set of values of $\delta_h^{ L_k }$ for all the combinations of selected remainders in the sequences that form the partial sum $S_h$ (see~\cite[Lemma 7.3]{rgb}); in the same way, we use the notation $\{ \delta_h^{ R_k } \}$ to denote the set of values of $\delta_h^{ R_k }$.

% Lema 1.3

\begin{lema} \label{S1biyeccion}

There is a bijective function $f_h: \{ \delta_h^{ L_k } \} \to \{ \delta_h^{ R_k } \}$ such that

\begin{align*}
f_h \left ( x \right ) = \delta_h - \left ( x - \delta_h \right ) \frac{ p_k^2 }{ m_k - p_k^2 },
\end{align*}

and

\begin{align*}
f_h^{-1} \left ( x \right ) = \delta_h + \left ( \delta_h - x \right ) \frac{ m_k - p_k^2 }{ p_k^2 }.
\end{align*}

\begin{proof}

For a given level $h \; (1 \le h \le k)$, if we change the combination of selected remainders in the partial sum $S_h$, some permitted $h$-tuples will be transferred from the Left interval $I[1, p_k^2]_h$ to the Right interval $I[p_k^2+1, m_k]_h$, or vice versa, as we have seen in the paragraph below~\cite[Definition 6.4]{rgb}. However, there exists a one-to-one correspondence between the size of the set of permitted $h$-tuples within the Left interval $I[1, p_k^2]_h$ and the size of the set of permitted $h$-tuples within the Right interval $I[p_k^2+1, m_k]_h$, since the size of the set of permitted $h$-tuples within the interval $I[1, m_k]_h$ is the same, regardless of the combination of selected remainders in the sequences $s_h$ that form the partial sum $S_h$, by~\cite[Proposition 2.3]{rgb} and~\cite[Lemma 6.1]{rgb}. It follows that there is also a one-to-one correspondence between the set of values of $\delta_h^{ L_k }$ and the set of values of $\delta_h^{ R_k }$. Therefore, for a given level $h \; (1 \le h \le k)$, we can define a bijective function $f_h: \{ \delta_h^{ L_k } \} \to \{ \delta_h^{ R_k } \}$.

Now, for the partial sum $S_h$, assume that the $h$-density within both $I [1, p_k^2 ]_h$ and $I [p_k^2+1, m_k ]_h$ is equal to the average $\delta_h$. Then, suppose that some permitted $h$-tuples are transferred from the Right interval to the Left interval. We have an increase $( \delta_h^{ L_k } - \delta_h )$ in the $h$-density within the Left interval and a decrease $( \delta_h - \delta_h^{ R_k } )$ in the $h$-density within the Right interval. See~\cite[(21)]{rgb}. Since there are $p_k^2 / p_h$ subintervals of size $p_h$ within the Left interval $I [1, p_k^2 ]_h$, by definition, the number of permitted $h$-tuples entering the Left interval is equal to $( \delta_h^{ L_k } - \delta_h ) p_k^2 / p_h$. In the same way, within the Right interval $I[p_k^2+1, m_k ]_h$ there are $( m_k - p_k^2 ) / p_h$ subintervals of size $p_h$, so, the number of permitted $h$-tuples exiting the Right interval is $( \delta_h - \delta_h^{ R_k } ) ( m_k - p_k^2 ) / p_h$. Since the number of permitted $h$-tuples entering the Left interval must be equal to the number of permitted $h$-tuples exiting the Right interval,

\begin{align*}
\left ( \delta_h^{ L_k } - \delta_h \right ) \frac{ p_k^2 }{ p_h } = \left ( \delta_h - \delta_h^{ R_k } \right ) \frac{ m_k - p_k^2 }{ p_h } \Longrightarrow \delta_h^{ R_k } = \delta_h - \left ( \delta_h^{ L_k } - \delta_h \right ) \frac{ p_k^2 }{ m_k - p_k^2 }.
\end{align*}

Therefore, we have a bijective function $f_h: \{ \delta_h^{ L_k } \} \to \{ \delta_h^{ R_k } \}$, such that $f_h ( x ) = \delta_h - ( x - \delta_h ) ( p_k^2 / ( m_k - p_k^2 ) )$, and we can check that $f_h^{-1} ( x ) = \delta_h + ( \delta_h - x ) ( m_k - p_k^2 ) / p_k^2$.

\end{proof}
\end{lema}

\begin{nota} \label{S1nota1}

Note that for a given level $h \; (1 \le h \le k)$, the image of $\max \{ \delta_h^{ L_k } \}$ under the function $f_h$ of the preceding lemma is $\min \{ \delta_h^{ R_k } \}$, and the image of $\min \{ \delta_h^{ L_k } \}$ under $f_h$ is $\max \{ \delta_h^{ R_k } \}$. See the proof of Lemma~\ref{S1biyeccion} and~\cite[(21)]{rgb}.

\end{nota}

By~\cite[Theorem 4.3]{rgb}, the average density of permitted $h$-tuples within every Left interval $I[1, p_k^2]_h$ is equal to $\delta_h$, that is, is equal to the $h$-density within the period of the partial sum $S_h$. The following lemma shows that the values of the $h$-density within the Left interval $I[1, p_k^2]_h$ of every partial sum $S_h \; (1 \le h \le k)$ converge uniformly to the respective average $\delta_h$ as $k \to \infty$, no matter the combination of selected remainders in every partial sum $S_k$.

\begin{nota} \label{S1nota2}

As we go from $h=1$ to $h=k$, the number of combinations of selected remainders within the periods of all the sequences $s_h \; (1 \le h \le k)$ forming $S_k$, given by

\begin{align*}
\displaystyle{p_1 \choose 1} \displaystyle{p_2 \choose 2} \displaystyle{p_3 \choose 2} \cdots \displaystyle{p_k \choose 2},
\end{align*}

increases (see~\cite[Definition 4.1]{rgb}); furthermore, the number of permitted $h$-tuples within the Right interval $I[p_k^2+1, m_k]_h$ decreases. Thus, it appears reasonable to assume that the variation of $\delta_h^{ R_k }$ about the mean $\delta_h$ increases from $h=1$ to $h=k$.

\end{nota}

% Definicion 1.1

\begin{definicion} \label{S1defin1}

We define the extended function $\varphi_h : \mathbb{R}_+ \to \mathbb{R}_+$ such that $\varphi_h = f_h$, where $f_h$ is the first function of Lemma~\ref{S1biyeccion}, and we define the extended function $\varphi_h^{-1} : \mathbb{R}_+ \to \mathbb{R}_+$ such that $\varphi_h^{-1} = f_h^{-1}$, where $f_h^{-1}$ is the second function of Lemma~\ref{S1biyeccion}, for every $h \; (1 \le h \le k)$.

\end{definicion}

\begin{nota} \label{S1nota3}

Note that for a given level $h \; (1 \le h \le k)$, the image of an upper bound for $\delta_h^{ L_k }$ under the function $\varphi_h$ is a lower bound for $\delta_h^{ R_k }$, and the image of a lower bound for $\delta_h^{ L_k }$ under $\varphi_h$ is an upper bound for $\delta_h^{ R_k }$. See the proof of Lemma~\ref{S1biyeccion} and~\cite[(21)]{rgb}.

\end{nota}

% Lema 1.4

\begin{lema} \label{S1deltahLkTodeltah}

Let $S_k \; (k \ge 4)$ be a partial sum of the series $\sum s_k$. Let us consider the Left interval $I[1, p_k^2]_h$ in every partial sum $S_h$ from $h = 1$ to $h = k$. For every small number $\Delta > 0$, there exists $N$ (depending only on $\Delta$) such that $k > N$ implies $| \delta_h^{ L_k } - \delta_h | < \Delta$ for every partial sum $S_h$ from $h = 1$ to $h = k$, regardless of the combination of selected remainders in $S_k$.

\begin{proof}

\begin{enumerate}[Step 1.]

\item Let us choose $h=h^\prime \; (1 < h^\prime < k)$ as a fixed level and let $\Delta > 0$ be a small number. As $k \to \infty$, for every level from $h = 1$ to $h = h^\prime$ the size of the Left interval $I[1, p_k^2]_h$ increases, so, the values of $\delta_h^{ L_k }$ converge uniformly to the respective average $\delta_h$, by~\cite[Proposition 5.3]{rgb}. Therefore, there exists $K$ such that

\begin{align}  \label{S1EQformu01}
\delta_h - \Delta < \min \{ \delta_h^{ L_k } \} \le \delta_h^{ L_k } \le \max \{ \delta_h^{ L_k } \} < \delta_h + \Delta
\end{align}

when $k > K$, for every Left interval $I[1, p_k^2]_h$ from $h = 1$ to $h = h^\prime$, regardless of the combination of selected remainders in $S_k$.

\item Let $\mathscr{P}$ be the sequence of primes. Given the small number $\Delta$, we define the sequence $\{ \epsilon_k \} = \epsilon_1, \epsilon_2, \epsilon_3, \epsilon_4, \ldots$ by the equations

\begin{align*}
&\epsilon_1 = 1,\\
&\epsilon_2= 1,\\
&\epsilon_3 = 1,\\
\text{and}\\
&\epsilon_k = \Delta \frac{p_k^2}{m_k - p_k^2} \qquad \text{for } k \ge 4, p_k \in \mathscr{P}.
\end{align*}

Note that $\epsilon_k$ decreases and tends to $0$ as $k \to \infty$.

Now, let us consider the Right interval $I[p_k^2+1, m_k]_h$ in every partial sum $S_h$ from level $h = 1$ to level $h = h^\prime$. Applying the function $\varphi_h$ to each member of \eqref{S1EQformu01} and using remarks~\ref{S1nota1} and~\ref{S1nota3}, we obtain

\begin{align*}
\varphi_h \left ( \delta_h - \Delta \right ) > \varphi_h \left ( \min \{ \delta_h^{ L_k } \} \right ) \ge \varphi_h \left ( \delta_h^{ L_k } \right ) \ge \varphi_h \left ( \max \{ \delta_h^{ L_k } \} \right ) > \varphi_h \left ( \delta_h + \Delta \right )
\end{align*}

if $k > K$, for every partial sum $S_h$ from $h = 1$ to $h = h^\prime$, regardless of the combination of selected remainders in $S_k$. Then, since the function $f_h$ is the restriction of $\varphi_h$ to the set $\{ \delta_h^{ L_k } \}$, we can write

\begin{align*}
\varphi_h \left ( \delta_h - \Delta \right ) > f_h \left ( \min \{ \delta_h^{ L_k } \} \right ) \ge f_h \left ( \delta_h^{ L_k } \right ) \ge f_h \left ( \max \{ \delta_h^{ L_k } \} \right ) > \varphi_h \left ( \delta_h + \Delta \right ),
\end{align*}

so, using the sequence $\{ \epsilon_k \}$ we obtain

\begin{align} \label{S1EQformu02}
\delta_h + \epsilon_k > \max \{ \delta_h^{ R_k } \} \ge \delta_h^{ R_k } \ge \min \{ \delta_h^{ R_k } \} > \delta_h - \epsilon_k
\end{align}

for the Right interval $I[p_k^2+1, m_k]_h$ in every partial sum $S_h$ from level $h = 1$ to level $h = h^\prime$ and $k > K$, regardless of the combination of selected remainders in $S_k$.

Note that given $\Delta$ (for the Left interval $I[1, p_k^2]_h$), the bounds $\delta_h - \epsilon_k$ and $\delta_h + \epsilon_k$ (for the Right interval $I[p_k^2+1, m_k]_h$) approximates $\delta_h$ more and more closely as $k$ increases, for every level from $h = 1$ to $h = h^\prime$, since $\{ \epsilon_k \}$ is a decreasing sequence.

\item Now, given a level $h \; (1 \le h \le h^\prime)$, let us consider the quotient $(\max \{ \delta_h^{ R_k } \} - \delta_h) / \epsilon_k$ for the Right interval $I[p_k^2+1, m_k]_h$. Using the bijective function $f_h$ we can write

\begin{align*}
\frac{\max \{ \delta_h^{ R_k } \} - \delta_h}{\epsilon_k} = \frac{f_h \left ( \min \{ \delta_h^{ L_k } \} \right ) - \delta_h}{\epsilon_k} = \frac{ \left ( \delta_h - \min \{ \delta_h^{ L_k } \} \right ) \frac{p_k^2}{m_k - p_k^2} }{ \Delta \frac{p_k^2}{m_k - p_k^2} },
\end{align*}

thus, we obtain

\begin{align}  \label{S1EQformu03}
\frac{\max \{ \delta_h^{ R_k } \} - \delta_h}{\epsilon_k} = \frac{ \left ( \delta_h - \min \{ \delta_h^{ L_k } \} \right ) }{ \Delta }.
\end{align}

In the same way it can be checked that

\begin{align}  \label{S1EQformu04}
\frac{ \delta_h - \min \{ \delta_h^{ R_k } \} }{\epsilon_k} = \frac{ \left ( \max \{ \delta_h^{ L_k } \} - \delta_h \right ) }{ \Delta }.
\end{align}

By~\cite[Proposition 5.3]{rgb} (see the explanation given in Step 1), the right-hand side of \eqref{S1EQformu03} and \eqref{S1EQformu04} tend to $0$ as $k \to \infty$, so, it follows that $(\max \{ \delta_h^{ R_k } \} - \delta_h)$ and $(\delta_h - \min \{ \delta_h^{ R_k } \})$ are both ${\rm o } ( \epsilon_k )$ as $k \to \infty$.

\item By Step 2, the inequalities in \eqref{S1EQformu02} are satisfied for the Right interval $I[p_k^2+1, m_k]_h$ in every partial sum $S_h$ from level $h = 1$ to level $h = h^\prime$, for $k > K$. In this and the next step we prove that these inequalities are also satisfied from level $h = h^\prime + 1$ to level $h = k$, for sufficiently large $k$.

We make the following remark: At the same time that the values of $\max \{ \delta_h^{ R_k } \}$ and $\min \{ \delta_h^{ R_k } \}$ from $h = 1$ to $h = h^\prime$ converge uniformly to $\delta_h$ as $k \to \infty$, by Lemma~\ref{S1deltahRkTodeltah}, the values of $\delta_h + \epsilon_k$ and $\delta_h - \epsilon_k$ tend also to $\delta_h$ for every level from $h = 1$ to $h = h^\prime$, since $\epsilon_k \to 0$ as $k \to \infty$, as we have seen in the last paragraph of Step 2. However, $\max \{ \delta_h^{ R_k } \}$ and $\min \{ \delta_h^{ R_k } \}$ converge to $\delta_h$ more rapidly than $\delta_h + \epsilon_k$ and $\delta_h - \epsilon_k$ respectively, as $k \to \infty$, for every level from $h = 1$ to $h = h^\prime$, by the result of the preceding step.

\item Now, from $h = 1$ to $h = k$ the values of $\max \{ \delta_h^{ R_k } \}$ and $\min \{ \delta_h^{ R_k } \}$ approximate $\delta_h$ more and more closely as $k$ increases, by Lemma~\ref{S1deltahRkTodeltah}. On the other hand, $\max \{ \delta_h^{ R_k } \}$ and $\min \{ \delta_h^{ R_k } \}$ converge to $\delta_h$ faster than $\delta_h + \epsilon_k$ and $\delta_h - \epsilon_k$ respectively, as $k \to \infty$, for every level from $h = 1$ to $h = h^\prime$, by the remark in the preceding step. Therefore, for every sufficiently large $k$ there must exist some $\epsilon \; (0 < \epsilon < \epsilon_k)$ such that $\max \{ \delta_h^{ R_k } \} < \delta_h + \epsilon < \delta_h + \epsilon_k$ and $\min \{ \delta_h^{ R_k } \} > \delta_h - \epsilon > \delta_h - \epsilon_k$ for every level from $h = 1$ to $h = k$, since the values of $\max \{ \delta_h^{ R_k } \}$ and $\min \{ \delta_h^{ R_k } \}$ converge uniformly to $\delta_h$ from $h = 1$ to $h = k$, by Lemma~\ref{S1deltahRkTodeltah}. Then there exists $K^\prime > K$ such that $\max \{ \delta_h^{ R_k } \} < \delta_h + \epsilon_k$ and $\min \{ \delta_h^{ R_k } \} > \delta_h - \epsilon_k$ for every level from $h = 1$ to $h = k$, for every $k > K^\prime$. Thus, we can write

\begin{align} \label{S1EQformu05}
\delta_h + \epsilon_k > \max \{ \delta_h^{ R_k } \} \ge \delta_h^{ R_k } \ge \min \{ \delta_h^{ R_k } \} > \delta_h - \epsilon_k
\end{align}

for the Right interval $I[p_k^2+1, m_k]_h$ in every partial sum $S_h$ from level $h = 1$ to level $h = k$, when $k > K^\prime$, regardless of the combination of selected remainders in $S_k$.

\item Now, applying the function $\varphi_h^{-1}$ to each member of \eqref{S1EQformu05} and using remarks~\ref{S1nota1} and~\ref{S1nota3}, we obtain

\begin{align} \label{S1EQformu06}
\varphi_h^{-1} \left ( \delta_h + \epsilon_k \right ) < \varphi_h^{-1} \left ( \max \{ \delta_h^{ R_k } \} \right ) \le \varphi_h^{-1} \left ( \delta_h^{ R_k } \right ) \le \varphi_h^{-1} \left ( \min \{ \delta_h^{ R_k } \} \right ) < \varphi_h^{-1} \left ( \delta_h - \epsilon_k \right )  \quad (k > K^\prime),
\end{align}

for every Left interval $I[1, p_k^2]_h$ from level $h = 1$ to level $h = k$, regardless of the combination of selected remainders in $S_k$. Now, since the function $f_h^{-1}$ is the restriction of $\varphi_h^{-1}$ to the set $\{ \delta_h^{ R_k } \}$, we can write

\begin{align} \label{S1EQformu07}
\varphi_h^{-1} \left ( \delta_h + \epsilon_k \right ) < f_h^{-1} \left ( \max \{ \delta_h^{ R_k } \} \right ) \le f_h^{-1} \left ( \delta_h^{ R_k } \right ) \le f_h^{-1} \left ( \min \{ \delta_h^{ R_k } \} \right ) < \varphi_h^{-1} \left ( \delta_h - \epsilon_k \right )  \quad (1 \le h \le k, k > K^\prime).
\end{align}

Then, by definition we have

\begin{align} \label{S1EQformu08}
\varphi_h^{-1} \left ( \delta_h + \Delta \frac{p_k^2}{m_k - p_k^2} \right ) < f_h^{-1} \left ( \max \{ \delta_h^{ R_k } \} \right ) \le f_h^{-1} \left ( \delta_h^{ R_k } \right ) \le f_h^{-1} \left ( \min \{ \delta_h^{ R_k } \} \right ) < \\ \nonumber < \varphi_h^{-1} \left ( \delta_h - \Delta \frac{p_k^2}{m_k - p_k^2} \right )  \quad (1 \le h \le k, k > K^\prime),
\end{align}

so,

\begin{align*}
\delta_h - \Delta < \min \{ \delta_h^{ L_k } \} \le \delta_h^{ L_k } \le \max \{ \delta_h^{ L_k } \} < \delta_h + \Delta  \quad (1 \le h \le k, k > K^\prime).
\end{align*}

From this it follows that

\begin{align*}
| \delta_h^{ L_k } - \delta_h | < \Delta,
\end{align*}

for the Left interval $I[1, p_k^2]_h$ in every partial sum $S_h$ from $h = 1$ to $h = k$, for every $k > N = K^\prime$, regardless of the combination of selected remainders in $S_k$. In words, while the $h$-density in every Right interval $I[p_k^2+1, m_k]_h$ from $h = 1$ to $h = k$ approximates the corresponding average $\delta_h$ more and more closely as $k \to \infty$ (by Lemma~\ref{S1deltahRkTodeltah}), the $h$-density in every Left interval $I[1, p_k^2+1]_h \; (1 \le h \le k)$ approximates $\delta_h$ as well. The lemma is proved.

\end{enumerate}

\end{proof}
\end{lema}

Note that using the preceding lemma and~\cite[Theorem 3.4]{rgb} the crucial result~\cite[Lemma 7.1]{rgb} can be easily proved.

% seccion 2
%-----------------------------------------------------------------------------------------------

\section{ A lower bound for $\delta_k^{L_k}$ } \label{seccion02}

Let $S_k \; (k \ge 4)$ be a partial sum of the series $\sum s_k$. Let us consider the Left interval $I[1, p_k^2]_h$ and the Right interval $I[p_k^2+1, m_k]_h$ in every partial sum $S_h$ from level $h = 1$ to level $h = k$.

Consider a fixed level $h^\prime \ge 4$. By~\cite[Lemma 3.2]{rgb}, the factor by which we must multiply the $h$-density within the interval $I[1, m_k]_h$ (denoted by $\delta_h$), to obtain the $(h+1)$-density within the interval $I[1, m_k]_{h+1}$ (denoted by $\delta_{h+1}$), for every level transition $h \to h+1$ from $h=1$ to $h=k-1$ is $(p_{h+1} -2) / p_h$, and this is the `average' case. Using this factor at each level transition from $h = h^\prime$ to $h = k-1$, we can write

\begin{align} \label{S2EQformu01}
\delta_k = \delta_{h^\prime} \left ( \frac{p_{h^\prime+1} -2}{p_{h^\prime}} \right ) \left ( \frac{p_{h^\prime+2} -2}{p_{h^\prime+1}} \right ) \cdots \left ( \frac{p_k -2}{p_{k-1}} \right ) = \delta_{h^\prime} \displaystyle\prod_{h=h^\prime}^{k-1} \left ( \frac{p_{h+1} -2}{p_h} \right ).
\end{align}

Now, suppose given an specific combination of selected remainders in the sequences $s_h \; (1 \le h \le k)$ that form $S_k$. Let $\delta_h^{ L_k }$, $\delta_h^{ R_k }$ be the $h$-density of permitted $h$-tuples in every interval $I[1, p_k^2]_h$, $I[p_k^2+1, m_k]_h$ respectively.

Let us denote by $\phi_h^{ L_k }$ ($\phi_h^{ R_k }$) the `true' factor by which we must multiply the $h$-density within $I[1, p_k^2]_h$ ($I[p_k^2+1, m_k]_h$) to obtain the $(h+1)$-density within $I[1, p_k^2]_{h+1}$ ($I[p_k^2+1, m_k]_{h+1}$), for every level transition $h \to h+1$ from $h=1$ to $h=k-1$. In symbols

\begin{align*}
&\delta_{h+1}^{ L_k } = \delta_h^{ L_k } \phi_h^{ L_k }, \\
&\delta_{h+1}^{ R_k } = \delta_h^{ R_k } \phi_h^{ R_k } \quad (1 \le h <k).
\end{align*}

Therefore, in the case of the Left intervals $I[1, p_k^2]_h$, we can write

\begin{align*}
\delta_k^{ L_k } = \delta_{h^\prime}^{ L_k } \phi_{h^\prime}^{ L_k } \phi_{h^\prime+1}^{ L_k } \cdots \phi_{k-1}^{ L_k } = \delta_{h^\prime}^{ L_k } \displaystyle\prod_{h={h^\prime}}^{k-1} \phi_h^{ L_k },
\end{align*}

and in the case of the Right intervals $I[p_k^2+1, m_k]_h \; (1 \le h \le k)$, we can write

\begin{align*}
\delta_k^{ R_k } = \delta_{h^\prime}^{ R_k } \phi_{h^\prime}^{ R_k } \phi_{h^\prime+1}^{ R_k } \cdots \phi_{k-1}^{ R_k } = \delta_{h^\prime}^{ R_k } \displaystyle\prod_{h={h^\prime}}^{k-1} \phi_h^{ R_k }.
\end{align*}

(Compare this formulas to \eqref{S2EQformu01}, given for the `average' case.)

Note that for every level transition $h \to h+1 \; (1 \le h < k)$, the factor $\phi_h^{ L_k }$ ($\phi_h^{ R_k }$) can be greater or lesser than the corresponding average factor $(p_{h+1} -2) / p_h$ given by~\cite[Lemma 3.2]{rgb} (see~\cite[Remark 7.2]{rgb}). Furthermore, it is obvious that the set of factors $\phi_{h^\prime}^{ L_k }, \phi_{h^\prime+1}^{ L_k }, \phi_{h^\prime+2}^{ L_k }, \ldots, \phi_{k-1}^{ L_k }$ ($\phi_{h^\prime}^{ R_k }, \phi_{h^\prime+1}^{ R_k }, \phi_{h^\prime+2}^{ R_k }, \ldots, \phi_{k-1}^{ R_k }$) depends on the combination of selected remainders in $S_k$.

\begin{nota} \label{S2nota1}

Let us consider the Left interval $I[1, p_k^2]_h$ of every partial sum $S_h$ from $h=1$ to $h=k$, where $k$ is sufficiently large. Since the selected remainders in every sequence $s_h$ remove permitted $(h-1)$-tuples from the partial sum $S_{h-1}$, the number of permitted $h$-tuples in $I[1, p_k^2]_h$ decreases from level $h=1$ to level $h=k$. From this, one might wrongly think that the rate of increase of the $h$-density in $I[1, p_k^2]_h$ between $h=1$ and $h=k$ could be less than the rate of increase of the average $\delta_h$, for a given combination of selected remainders in $S_k$. In other words, suppose we assume the following: There exists a particular combination of selected remainders in the sequences $s_h \; (1 \le h \le k )$ that form $S_k$ such that $\delta_h^{ L_k } < \delta_h$ for every level from $h=1$ to $h=k$, for sufficiently large $k$. Clearly the preceding assumption contradicts~\cite[Remark 7.2]{rgb}. (The same can be said with regard to the Right interval $I[p_k^2+1, m_k]_h$ of every partial sum $S_h$ from $h=1$ to $h=k$, where $k$ is sufficiently large.)

\end{nota}

Now, in the case of the Right intervals $I[p_k^2+1, m_k]_h  \; (1 \le h \le k)$ we might (wrongly) assume that for every $k$ there exists a combination of selected remainders in the sequences that form $S_k$ and a corresponding set of factors $\phi_{h^\prime}^{ R_k }, \phi_{h^\prime+1}^{ R_k }, \phi_{h^\prime+2}^{ R_k }, \ldots, \phi_{k-1}^{ R_k }$, such that

\begin{align*}
\max \{ \delta_k^{ R_k } \} = \max \{ \delta_{h^\prime}^{ R_k } \} \phi_{h^\prime}^{ R_k } \phi_{h^\prime+1}^{ R_k } \cdots \phi_{k-1}^{ R_k }.
\end{align*}

Clearly, for sufficiently large $k$ this is not possible, by the preceding remark, that is, the combination of selected remainders in $S_k$ for the maximum value of the $h^\prime$-density in the Right interval $I[p_k^2+1, m_k]_{h^\prime}$ is not necessarily the same as the combination for the maximum value of the $k$-density in $I[p_k^2+1, m_k]_k$. (A similar formula for $\min \{ \delta_k^{ R_k } \}$ is not possible either, for the same reasons.)

However, for sufficiently large $k$ the values of $\delta_h^{ R_k } \; (1 \le h \le k)$ are very close to the average $\delta_h$ by Lemma~\ref{S1deltahRkTodeltah}, so, for every level transition $h \to h+1 \; (1 \le h < k)$ the factor $\phi_h^{ R_k }$ will be very close to the average factor $(p_{h+1} -2) / p_h$, regardless of the combination of selected remainders in $S_k$. Therefore, it seems reasonable to think the following: Between $h=h^\prime$ and $h=k$, the increase from $\max \{ \delta_{h^\prime}^{ R_k } \}$ to $\max \{ \delta_k^{ R_k } \}$ (and the increase from $\min \{ \delta_{h^\prime}^{ R_k } \}$ to $\min \{ \delta_k^{ R_k } \}$) is approximately proportional to the increase from $\delta_{h^\prime}$ to $\delta_k$.

Likewise, in the case of the Left intervals $I[1, p_k^2]_h  \; (1 \le h \le k)$, we can say that between $h=h^\prime$ and $h=k$, the increase from $\max \{ \delta_{h^\prime}^{ L_k } \}$ to $\max \{ \delta_k^{ L_k } \}$ (and the increase from $\min \{ \delta_{h^\prime}^{ L_k } \}$ to $\min \{ \delta_k^{ L_k } \}$) is approximately proportional to the increase from $\delta_{h^\prime}$ to $\delta_k$, since for sufficiently large $k$ the values of $\delta_h^{ L_k } \; (1 \le h \le k)$ are very close to the average $\delta_h$, by Lemma~\ref{S1deltahLkTodeltah}.

These facts allows us to establish a lower bound for $\delta_k^{ L_k }$ in the following lemma.

% Lema 2.1

\begin{lema} \label{S2FormuPropor}

Let $S_k \; (k > 4)$ be a partial sum of the series $\sum s_k$. Let us consider the Left interval $I[1, p_k^2]_h$ in every partial sum $S_h$ from $h = 1$ to $h = k$. Given a fixed level $h=h^\prime \; (4 \le h^\prime < k)$ and a small real number $\Delta > 0$, we have

\begin{align} \label{S2EQformu02}
\delta_k^{ L_k } > \delta_k - \frac{\delta_k}{\delta_{h^\prime}} \Delta
\end{align}

for sufficiently large $k$, regardless of the combination of selected remainders in $S_k$.

\begin{proof}

Let us consider the Left interval $I[1, p_k^2]_h$ in every partial sum $S_h$ from level $h = 1$ to level $h = k$. By Lemma~\ref{S1deltahLkTodeltah}, given a small number $\Delta > 0$ there exists $N$ (depending only on $\Delta$) such that

\begin{align*}
\delta_h - \Delta < \delta_h^{ L_k } < \delta_h + \Delta \quad (1 \le h \le k)
\end{align*}

for every $k > N$, regardless of the combination of selected remainders in $S_k$. Since $\delta_k / \delta_{h^\prime} > 1$, it follows that

\begin{align*}
\delta_k^{ L_k } > \delta_k - \Delta > \delta_k - \frac{\delta_k}{\delta_{h^\prime}} \Delta  \quad (k > N)
\end{align*}

and we are done.

\end{proof}
\end{lema}

\begin{proof}[Alternative proof of Lemma 2.1]

\begin{enumerate}[Step 1.]

\item Let us consider the Right interval $I[p_k^2+1, m_k]_h$ in every partial sum $S_h$ from level $h = 1$ to level $h = k$. By the Step 4 in Lemma~\ref{S1deltahLkTodeltah}, given a small number $\Delta > 0$ and the associated sequence $\{ \epsilon_k \}$, there exists $N$ such that

\begin{align} \label{S2EQformu03}
\delta_h - \epsilon_k < \delta_h^{ R_k } < \delta_h + \epsilon_k \quad (1 \le h \le k)
\end{align}

when $k > N$, regardless of the combination of selected remainders in $S_k$. Now, from \eqref{S2EQformu03} it follows that

\begin{align*}
\delta_k^{ R_k } < \delta_k + \epsilon_k < \delta_k + \frac{\delta_k}{\delta_{h^\prime}} \epsilon_k  \quad (k > N),
\end{align*}

since $\delta_k / \delta_{h^\prime} > 1$. Then, by the definition of the sequence $\{ \epsilon_k \}$ we obtain the upper estimate

\begin{align}  \label{S2EQformu04}
\delta_k^{ R_k } < \delta_k + \frac{\delta_k}{\delta_{h^\prime}} \Delta \frac{p_k^2}{m_k - p_k^2}  \quad (k > N).
\end{align}

\item Now, applying the function $\varphi_k^{-1}$, to each member of \eqref{S2EQformu04} and using Remark~\ref{S1nota3}, we obtain

\begin{align*}
\varphi_k^{-1} \left ( \delta_k^{ R_k } \right ) > \varphi_k^{-1} \left ( \delta_k + \frac{\delta_k}{\delta_{h^\prime}} \Delta \frac{p_k^2}{m_k - p_k^2} \right )  \quad (k > N),
\end{align*}

regardless of the combination of selected remainders in $S_k$. Then, since the function $f_k^{-1}$ is the restriction of $\varphi_k^{-1}$ to the set $\{ \delta_k^{ R_k } \}$, we can write

\begin{align*}
f_k^{-1} \left ( \delta_k^{ R_k } \right ) > \varphi_k^{-1} \left ( \delta_k + \frac{\delta_k}{\delta_{h^\prime}} \Delta \frac{p_k^2}{m_k - p_k^2} \right )  \quad (k > N),
\end{align*}

thus, in the Left interval $I[1, p_k^2]$ we obtain the lower estimation

\begin{align*}
\delta_k^{ L_k } > \delta_k - \frac{\delta_k}{\delta_{h^\prime}} \Delta  \quad (k > N).
\end{align*}

\end{enumerate}

\end{proof}

\begin{figure}[H]
\centering
\includegraphics[width=10cm]{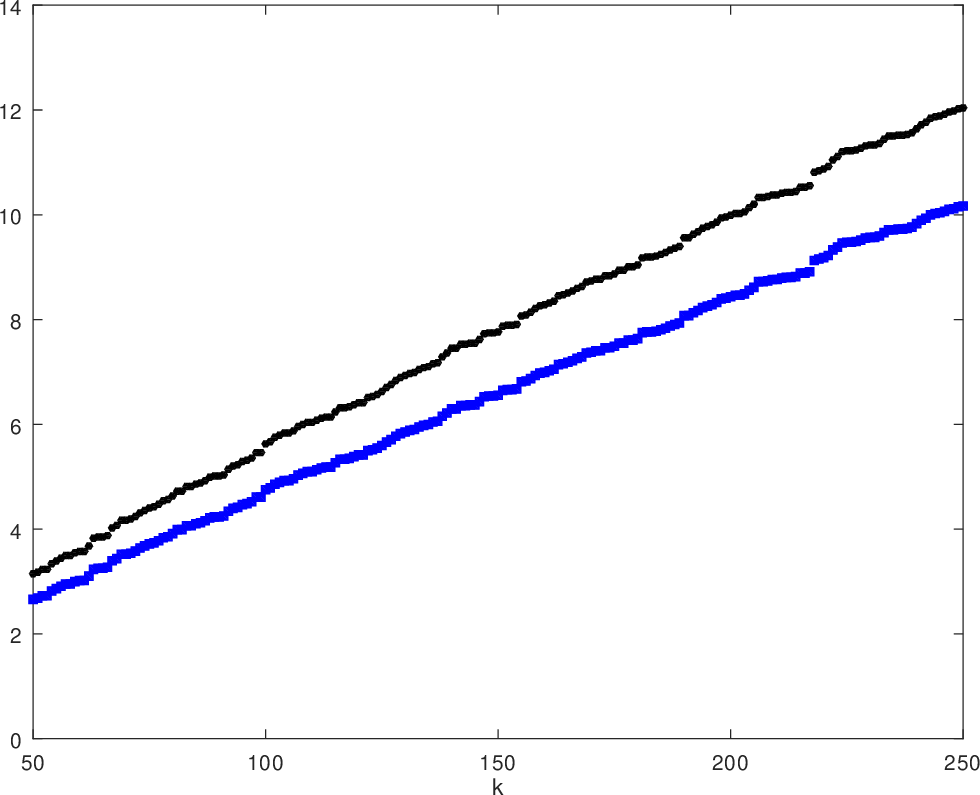}
\caption{$\delta_k$ vs. right-hand side in \eqref{S2EQformu02} for $k$ ranging from 50 to 250.}
\label{Figura4}
\end{figure}

Taking $\Delta = 0.1$ and $h^\prime=5$ as a fixed level we can check, using~\cite[Lemma 5.2]{rgb}, that $\delta_h - \Delta < \delta_h^{L_k} < \delta_h + \Delta \; (1 \le h \le h^\prime)$, for every $k \ge 32$. Assume that $\delta_h - \Delta < \delta_h^{L_k} < \delta_h + \Delta$ is also true from $h=h^\prime+1$ to $h=k$ for $k \ge 32$. Thus, using this value of $\Delta$ in formula \eqref{S2EQformu02} we obtain the graphs of Figure \ref{Figura4}, where for each $k$ between $50$ and $250$ (the horizontal axis), a black dot represents $\delta_k$ and a blue box represents the right-hand side in \eqref{S2EQformu02}.

% seccion 3
%---------------------------------------------------------------------------------------------

\section{Estimating a value for $K_\alpha$} \label{seccion03}

Let us consider the partial sum $S_k \; (k \ge 4)$. Recall that~\cite[Theorem 7.1]{rgb} and~\cite[Corollary 7.2]{rgb} ensure the existence of a number $K_\alpha$ such that $\delta_k^{ L_k } > \delta_4$ for every $k > K_\alpha$, regardless of the combination of selected remainders in $S_k$.

Given a fixed level $h=h^\prime \; (1 \le h \le k)$, the $h$-density in the Left interval $I[1, p_k^2]_h$ (denoted by $\delta_h^{L_k}$) converges to the average $\delta_h$ for every level from $h=1$ to $h=h^\prime$ as $k \to \infty$, by~\cite[Proposition 5.3]{rgb}. In addition, given a small number $\Delta$, the~\cite[Proposition 5.2]{rgb} allows us to find $K$ such that for every $k > K$ we have $\delta_h - \Delta < \delta_h^{L_k} < \delta_h + \Delta$ for every level from $h=1$ to $h=h^\prime$ (see~\cite[Lemma 7.1, steps 3 and 4]{rgb}). However, these propositions do not guarantee that $\delta_h - \Delta < \delta_h^{L_k} < \delta_h + \Delta$ for every level from $h=1$ to $h=k \; (k > K)$, which is necessary for our purposes.

On the other hand, by Lemma~\ref{S1deltahLkTodeltah} we have $\delta_h - \Delta < \delta_h^{L_k} < \delta_h + \Delta$ for $h=1$ to $h=k$ and sufficiently large $k$, regardless of the combination of selected remainders in $S_k$. Furthermore, Lemma~\ref{S2FormuPropor} also establish a lower bound for $\delta_k^{L_k}$, for sufficiently large $k$. Now, the million dollar question is: How large should be the level $k$ for being considered `sufficiently large' in these cases?

Given the first period of the partial sum $S_k$, we have control over the $h$-density within the intervals $I[1, m_k]_h$ from $h=1$ to $h=k$, for every $k$; furthermore, we have control over the $h$-density within the Right intervals $I[p_k^2+1, m_k]_h \; (1 \le h \le k)$, for sufficiently large $k$. However, we have no control over the $h$-density within the Left intervals $I[1, p_k^2]_h$ for all the levels between $h=1$ and $h=k$. Therefore, to give an answer to the preceding question, it is first necessary `to take a detour' through the Right block of the partition, using the fact that $p_k^2 = {\rm o } ( m_k )$ (see the paragraph below Lemma~\ref{S1pk2o(ck)}), to find $K$ such that we can assume that $\delta_h^{R_k}$ is very close to the average $\delta_h$ from $h=1$ to $h=k$, for every $k > K$. After that we return to the Left block of the partition for proving that $\delta_k^{L_k} > \delta_4$ for every $k > K$.

% Lema 3.1

\begin{lema}[Estimation for $K_\alpha$] \label{SPestimacion} 

Let $S_k$ be a partial sum of the series $\sum s_k$. Let $K_\alpha$ be the number whose existence is guaranteed by ~\cite[Corollary 7.2]{rgb}. Then, if $K_\alpha = 2800$, for every $k > K_\alpha$ we have $\delta_k^{L_k} > \delta_4$, regardless of the combination of selected remainders in $S_k$.

\begin{proof}

\begin{enumerate}[Step 1.]

\item Given the partial sum $S_k \; (k > 5)$, let us take $h^\prime = 5$ as a fixed level. Furthermore let $K = 2800$ and $\Delta = 0.1$. We begin by establishing bounds for $\delta_h^{L_k}$ in every Left interval $I[1, p_k^2]_h$ from $h=1$ to $h=h^\prime$. Given a level $k > K$ we have $p_k^2 > 644702881$, thus, for the Left interval $I[1, p_k^2]_h$, using~\cite[Lemma 5.2]{rgb} we can check that

\begin{align} \label{S3EQformu02}
\delta_h - \Delta < \delta_h^{L_k} < \delta_h + \Delta,
\end{align}

for every level from $h=1$ to $h=h^\prime$ and $k > K$, regardless of the combination of selected remainders in $S_k$.

\item Applying the function $\varphi_h$ to every member of \eqref{S3EQformu02} and using Remark~\ref{S1nota3}, for $k > K$ we have

\begin{align*}
\varphi_h \left ( \delta_h - \Delta \right ) > \varphi_h \left ( \delta_h^{L_k} \right ) > \varphi_h \left ( \delta_h + \Delta \right ),
\end{align*}

so,

\begin{align}  \label{S3EQformu03}
\delta_h + \epsilon_k > \delta_h^{R_k} > \delta_h - \epsilon_k
\end{align}

for every Right interval $I[p_k^2+1, m_k]_h$ from $h=1$ to $h=h^\prime$, where $\epsilon_k \; (k > K)$ is a member of the sequence defined in the Step 2 of Lemma~\ref{S1deltahLkTodeltah}.

\item For a partial sum $S_k$ where $k > K$ the ratio of the size of the Left interval $I[1, p_k^2]$ to the size of the interval $I[ 1, m_k]$ is less than $1.4 \times 10^{-301}$. Thus, the size of the Left interval $I[ 1, p_k^2]$ is completely negligible compared to the size of $I[ 1, m_k]$. Hence, for every level from $h = 1$ to $h = k \; (k > K)$, in every Right interval $I[p_k^2+1, m_k]_h$ the value of $\delta_h^{R_k}$ will be very close to $\delta_h$ (almost equal to $\delta_h$), for every combination of selected remainders in $S_k$ (see the paragraph below Lemma~\ref{S1pk2o(ck)}).

\item Now, $\delta_h^{ R_k }$ tends uniformily to $\delta_h$ for every level from $h = 1$ to $h = k$, by Lemma~\ref{S1deltahRkTodeltah}, and simultaneously $\delta_h^{ L_k }$ tends uniformily to $\delta_h$ for every level from $h = 1$ to $h = k$, by Lemma~\ref{S1deltahLkTodeltah}, as $k \to \infty$. Therefore, it seems reasonable to think as follows: If for sufficiently large $k$ is $\delta_h^{R_k}$ very close to $\delta_h$ for every level from $h = 1$ to $h = k$, then $\delta_h^{L_k}$ will be quite close to $\delta_h$ for every level from $h = 1$ to $h = k$, regardless of the combination of selected remainders in $S_k$.

Since $\delta_h^{R_k}$ is almost equal to $\delta_h$ for every level from $h = 1$ to $h = k$ when $k > K$, by Step 3, we can assume that $\max \{ \delta_h^{ R_k } \} < \delta_h + \epsilon_k$ and $\min \{ \delta_h^{ R_k } \} > \delta_h - \epsilon_k$ for every level from $h = 1$ to $h = k$ whenever $k > K$ (see Step 5 in Lemma~\ref{S1deltahLkTodeltah}), that is, \eqref{S3EQformu03} is also satisfied for every level from $h = 1$ to $h = k$ when $k > K$. Then, by applying to \eqref{S3EQformu03} the procedure in Step 6 of the proof of Lemma~\ref{S1deltahLkTodeltah} we obtain

\begin{align*}
\delta_h - \Delta < \delta_h^{L_k} < \delta_h + \Delta \quad (1 \le h \le k, k > K),
\end{align*}

regardless of the combination of selected remainders in $S_k$; so, we can see that the result of Lemma~\ref{S1deltahLkTodeltah} holds for every $k > K$. In particular for $h=k$, we have the lower estimate

\begin{align}  \label{S3EQformu04}
\delta_k^{L_k} > \delta_k - \Delta \quad (k > K).
\end{align}

\item On the other hand, from \eqref{S3EQformu04} it follows that

\begin{align}  \label{S3EQformu05}
\delta_k^{L_k} > \delta_k - \frac{\delta_k}{\delta_{h^\prime}} \Delta
\end{align}

for every $k > K$, by the proof of Lemma~\ref{S2FormuPropor}.

\item Now, note that $\delta_k > 102.62$ and the ratio $\delta_k / \delta_{h^\prime}$ is greater than $159.63$ for $k > K$. Thus, using \eqref{S3EQformu04} we obtain

\begin{align}  \label{S3EQformu06}
\delta_k^{L_k} > \delta_k - \Delta = 102.52
\end{align}

and using \eqref{S3EQformu05} we obtain

\begin{align}  \label{S3EQformu07}
\delta_k^{L_k} > \delta_k - \frac{\delta_k}{\delta_{h^\prime}} \Delta = 86.66
\end{align}

for every $k > K$, regardless of the combination of selected remainders in $S_k$. Clearly, by the lower estimates \eqref{S3EQformu06} and \eqref{S3EQformu07} we conclude that $\delta_k^{ L_k } > \delta_5 > \delta_4$ for every level $k > K$, regardless of the combination of selected remainders in $S_k$. Thus, taking $K_\alpha = K = 2800$ the lemma is proved.

\end{enumerate}

\end{proof}

\end{lema}

% seccion 4
%--------------------------------------------------------------------------------------------

\section{Conclusion} \label{seccion04}

\begin{proof}[\bf Proof of the Main Theorem]

By~\cite[Lemma 7.1]{rgb} and~\cite[Corollary 7.2]{rgb}, there exists $K_\alpha$ such that $\delta_k^{ L_k } > \delta_4$ if $k > K_\alpha$, for every combination of selected remainders in the sequences $s_h$ that form the partial sum $S_k$. Furthermore, by Lemma~\ref{SPestimacion}, taking $K_\alpha = 2800$ the preceding statement is satisfied. Therefore, by~\cite[Lemma 8.3]{rgb}, every even integer $x > p_k^2$ is the sum of two primes, where $k > K_\alpha$. That is, if $x$ is greater than $p_{2800}^2 = 644702881$ it can be expressed as the sum of two primes.

Now, it is a known fact that the strong Goldbach conjecture has already been verified for all even numbers up to $4 \times 10^{18}$~\cite{oli}. Therefore, we conclude that every even number $x > 4$ can be expressed as the sum of two primes; thus, the binary Goldbach conjecture is proved and we are done.

\end{proof}

% Agradecimientos

\begin{acknowledgements}

The author wants to thank Dra. Patricia Quattrini (Departamento de Matematicas, FCEyN, Universidad de Buenos Aires) for helpful conversations and Dr. Hendrik W. Lenstra (Universiteit Leiden, The Netherlands) for his extremely useful suggestions.

\end{acknowledgements}

Ricardo G. Barca\\
Universidad Tecnologica Nacional\\
Buenos Aires (Argentina)\\
{\it E-mail address:} rbarca@frba.utn.edu.ar


{\small \begin{thebibliography}{99}

\bibitem{oli} Tom\'as Oliveira e Silva, Siegfried Herzog and Silvio Pardi, {\it Empirical verification of the even Goldbach conjecture and computation of prime gaps up to $4 \cdot 10^{18}$}, Mathematics of Computation, vol. 83, no. 288, pp. 2033-2060, Jul. 2014.

\bibitem{rgb} Ricardo Barca, Every sufficiently large even number is the sum of two primes (preprint),
\url{https://hal.archives-ouvertes.fr/hal-02075531}.

\end{thebibliography}}
\end{document}